\documentclass[a4, 12pt]{article}
\usepackage[table]{xcolor}
\usepackage[margin=1in]{geometry}
\usepackage{authblk}
\usepackage[numbers]{natbib}
\usepackage{longtable}
\usepackage{tikz}
\usepackage{subcaption}
\usepackage{caption}
\usepackage{float}
\usetikzlibrary{arrows.meta, positioning}
\usepackage{amsmath}
\usepackage{amssymb}
\usepackage{enumitem}

\title{Optimizing rake-links independently of timetables in railway operations}
\author[1]{Sourav Dey}
\affil[1]{Graduate School of Frontier Sciences, University of Tokyo}

\begin{document}

\maketitle

\begin{abstract}
This study addresses optimal rake-link formation in large-scale timetabled rail operations by modeling the problem as a directed acyclic graph and solving it via the minimum path cover algorithm. It enables efficient rake-to-service assignment while minimizing fleet size. Crucially, it decouples rake-link optimization from the timetable planning process, allowing planners to evaluate feasible rake configurations independently. The model incorporates operational constraints such as deadhead limits, service balance, and slack allowances. Applied to real-world data from Indian Railways, the results reveal clustered Pareto fronts in the decision space, indicating robust and redundant solutions. The approach lays a foundation for resilient, adaptive rail management via digital twin systems.
\end{abstract}

\section{Introduction}

Rake-links serve as the backbone of any efficient railway operation, as they govern how a finite pool of physical train units (rakes) is assigned to a succession of scheduled services. A rake represents a complete train set, comprising locomotives and coaches, such as the (Electric Multiple Unit) EMU, that must be routed from one service to the next so that no service remains unassigned and no rake remains idle on a running line for extended periods. In traditional planning, these assignments are made concurrently with the design of the timetable, creating a complex interdependency. The timetable presupposes specific rake availabilities, and rake assignments must accommodate every leg of the timetabled network. For example, if a late evening service is added between station A and station B, the schedule planner must immediately verify that a rake is available at or near station A at the needed time. Similarly, once a rake completes its designated run, it cannot remain on a main track waiting for its next assignment because prolonged idle time on a running line can block other services or incur storage penalties. So, timetable design and rake-link decisions are tightly coupled: changing a departure-arrival slot requires rechecking the corresponding rake assignment, and altering a rake assignment may force a timetable adjustment.

This coupling, while guaranteeing feasibility under steady conditions, becomes a handicap when planning needs to be agile or when alternative timetables must be evaluated quickly. If planners wish to test a revised evening‐peak schedule with slightly adjusted departure times or route patterns they must re‐evaluate the entire rake‐link configuration, often leading to lengthy iterations. Moreover, ad hoc demands from concerned stakeholders, such as rolling stock department requiring a rake withdrawal for maintenance or operations control wanting to adjust dwell times at intermediate stations, cannot be incorporated seamlessly unless both the timetable and rake-link are rechecked together. In large networks with hundreds of services and dozens of rakes, manual or semi-automated methods struggle to keep pace, stalling critical decision processes, and sometimes forcing suboptimal compromises.

Ideally, from the perspective of a railway operating company, it would be beneficial if these two major processes are decoupled so that the timetable can be designed in coordination with commercial, maintenance, and infrastructure stakeholders using heuristics or optimization routines while having a fast, on-demand feasibility check for rake-linking to flag conflicts or shortages. Under such an approach, a planner could propose a candidate schedule and immediately run a 'rake-link audit' that confirms the existing fleet suffices or highlights where additional rakes or deadhead movements are required. This capability would not only accelerate comparative studies of alternate schedules, such as testing a later departure from a suburban terminus or inserting an express run alongside existing locals, but also allow planners to explore “what-if” scenarios, such as temporarily withdrawing a rake for overhaul, without rebuilding the timetable from scratch. 

In this article, we aim to achieve this decoupling and focus on the rake-link optimization problem, given the time table and network topology. We begin with a formal characterization of the rake‐linking problem under typical operational constraints, including minimum layover times, allowable deadhead distances and speeds, and headway safety buffers. We then describe an algorithmic framework that can rapidly re‐evaluate rake feasibility as the timetable evolves. Finally, we propose a modular approach in which timetable optimization and rake assignment become loosely coupled subroutines, thus facilitating real-time 'what if' analyses and supporting network planners in making informed trade-offs between fleet utilization, service coverage, and operational robustness.

The remainder of this paper is structured as follows- Section \ref{sec:background} presents the background and review of the literature. Section \ref{sec:methodology} starts with the fleet minimization problem and later proposes a way to test a multi-objective variant that considers additional objectives in addition to the fleet size. Section \ref{sec:results} discusses different experimental observations and results. In Section \ref{sec:conclusion} we conclude with the limitations of the proposed approach along with suggestions for improvements and future research directions.

\section{Background}
\label{sec:background}
rake-linking is particularly challenging due to several operational constraints: minimum turnaround times for cleaning and boarding, compatibility of rakes with station platform lengths or track directions, potential for rake reversal, and maintenance requirements over periodic cycles. Moreover, the possibility of inserting empty (deadhead) runs between services adds another layer of complexity, especially in dense suburban or intercity networks. It is fundamentally an assignment problem in the domain of operations research, and as such, it has been extensively studied both in academia and practice.
In \cite{cordeau2001simultaneous} the simultaneous assignment of locomotives and passenger cars to scheduled trains is modeled as a multicommodity network flow problem (MCNF), where each type of locomotive and car is treated as a flow through a time-space graph, while \cite{ziarati1999branch} addresses the locomotive assignment problem using a novel `branch first cut second' approach. The train set assignment problem is also addressed in \cite{li2016pragmatic} through particle swarm optimization (PSO), along with maintenance schedule constraints related to a railway network in China. In \cite{rangaraj2006rake} the authors have solved the rake-linking as a minimum cost flow problem (MCF) and implemented their solution in a portion of the Indian railway network. They also perturb the parameters to identify areas of improvement, such as how a slight modification in the timetable can save an additional rake. The authors in \cite{biswas2024indian} approach the same problem of train-set assignment in the suburban operations of another division of Indian Railways using mixed-integer linear programming (MILP). In \cite{hong2009pragmatic}, authors present a network-based, heuristic method for the train-set assignment problem, emphasizing practical feasibility and fast runtime. The solution is applied to a real-life network in KTX network of S.Korea, that produced feasible, near-optimal solutions suitable for daily operations planning. In \cite{ghoseiri2004multi} first the Pareto frontiers are found through the $\varepsilon $- constraint method, and then multiobjective optimization is performed on the frontiers using three types of distance measures. A systematic review on the topic can be found in \cite{caprara2007passenger} and \cite{marques2023literature}.

The present work builds upon this body of literature by formulating the rake-linking problem as a minimum path cover problem over a directed acyclic graph. The current work is closely related to the work in \cite{rangaraj2006rake} and also uses the real life data from the same suburban railway network as in \cite{biswas2024indian}. First a single objective problem is proposed where the objective is to minimize the required fleet size. However, solutions with minimal rake usage often form tightly coupled service chains with limited buffer time, leading to:

\begin{itemize}
\item \textbf{Operational fragility}: Delays in one service can cascade to others.
\item \textbf{Imbalance}: Some rakes may be overused while others remain idle.
\item \textbf{Inflexibility}: Rigid links may break under disruptions.
\end{itemize}

Keeping this in mind, we finally extend the concept to a multi-objective variant, where other practical objectives are considered.

\section{Methodology}
\label{sec:methodology}
First we approach the single objective of minimizing fleet size through minimum path cover on directed acyclic graphs as follows:
\subsection{Fleet Size Minimization}
\label{lab:fleet_size_minimization}
\subsubsection*{Link Feasibility Graph}

Let \(\mathcal{S} = \{1,2,\dots,N\}\) denote the set of scheduled services. We construct a directed acyclic graph (DAG)
\[
  \mathcal{G}_F \;=\; (\mathcal{S},\,\mathcal{E}),
\]
called the \emph{Link Feasibility Graph},as follows:
\begin{itemize}
  \item Each node \(i \in \mathcal{S}\) represents service \(i\).
  \item There is a directed edge \((i \to j)\in \mathcal{E}\) if and only if a single rake can operate the service \(i\) immediately followed by the service \(j\). Equivalently,
    \[
      d_i \;=\; o_j
      \quad\text{and}\quad
      t^{\mathrm{arr}}_i + \delta_{ij} \;\le\; t^{\mathrm{dep}}_j,
    \]
    where:
    \begin{itemize}
      \item \(d_i\) is the destination station of service \(i\).
      \item \(o_j\) is the origin station of service \(j\).
      \item \(t^{\mathrm{arr}}_i\) is the scheduled arrival time of service \(i\).
      \item \(t^{\mathrm{dep}}_j\) is the scheduled departure time of service \(j\).
      \item \(\delta_{ij}\) is the minimum turnaround time required between services $i$ and $j$.
    \end{itemize}
\end{itemize}

\noindent
Since \(t^{\mathrm{arr}}_i < t^{\mathrm{dep}}_j\) whenever \((i\to j)\) exists, no directed cycle can exist. Hence \(\mathcal{G}_F\) is a DAG.

\subsubsection*{Minimum Path Cover and Fleet Size}

A \emph{path cover} of a DAG is a collection of vertex‐disjoint directed paths that together include every node exactly once. In our context, each path corresponds to a sequence of services handled by a single rake. Hence the minimum number of rakes required to cover all services is exactly the size of a \emph{Minimum Path Cover} \(MPC\)of \(\mathcal{G}_F\).

\[
  \boxed{
    \text{Fleet size}_{\min}
    \;=\; \bigl|\text{Minimum Path Cover}(\mathcal{G}_F)\bigr|
  }
\]

\subsubsection*{Reformulation via Bipartite Matching}

It is well known that finding a minimum path cover in a DAG of size \(\lvert \mathcal{S} \rvert\) is equivalent to computing a maximum matching in an associated bipartite graph, as follows:

\begin{enumerate}
  \item Create two copies of the vertex set:
    \[
      U \;=\; \{\,u_i : i \in \mathcal{S}\,\}, 
      \quad 
      V \;=\; \{\,v_i : i \in \mathcal{S}\,\}.
    \]
  \item For every directed edge \((i \to j)\in \mathcal{E}\), insert an undirected (bipartite) edge \((u_i,\,v_j)\).
  \item Call this bipartite graph \(\mathcal{B} = \bigl(U \cup V,\;\mathcal{E}'\bigr)\).
\end{enumerate}

Let \(\nu\) be the size of a maximum matching in \(\mathcal{B}\). Then the cardinality of a minimum path cover in the original DAG is:
\[
  \bigl|\text{Minimum Path Cover}(\mathcal{G}_F)\bigr|
  \;=\; |\mathcal{S}| \;-\;\nu.
\]
Consequently, the minimum number of rakes needed is
\[
  \min \lvert \text{Rakes} \rvert
  \;=\; |\mathcal{S}| \;-\;\nu,
\]
where \(\nu\) is computed via maximum bipartite matching.

\begin{figure}[H]
    \centering
    \begin{subfigure}[b]{0.23\textwidth}
        \centering
        \begin{tikzpicture}[->,>=stealth, node distance=1cm and 1cm]
            \tikzstyle{every node} = [draw, circle, minimum size=5mm, font=\scriptsize]
            
            \node (1) at (0,0) {1};
            \node (2) [right=of 1] {2};
            \node (3) [below=of 1] {3};
            \node (4) [right=of 3] {4};
            \node (5) [below=of 3] {5};
            \node (6) [right=of 5] {6};

            \draw (1) -- (2);
            \draw (1) -- (4);
            \draw (3) -- (4);
            \draw (3) -- (6);
            \draw (5) -- (6);
        \end{tikzpicture}
        \caption{\scriptsize Original DAG \(\mathcal{G}_F\)}
    \end{subfigure}
    \hfill
    \begin{subfigure}[b]{0.23\textwidth}
        \centering
        \begin{tikzpicture}[->,>=stealth, node distance=0.7cm and 0.5cm]
            \tikzstyle{every node} = [draw, circle, minimum size=5mm, font=\scriptsize]
            
            \node (u1) at (0,2.5) {$u_1$};
            \node (u2) at (0,1.5) {$u_2$};
            \node (u3) at (0,0.5) {$u_3$};
            \node (u4) at (0,-0.5) {$u_4$};
            \node (u5) at (0,-1.5) {$u_5$};
            \node (u6) at (0,-2.5) {$u_6$};

            \node (v1) at (2,2.5) {$v_1$};
            \node (v2) at (2,1.5) {$v_2$};
            \node (v3) at (2,0.5) {$v_3$};
            \node (v4) at (2,-0.5) {$v_4$};
            \node (v5) at (2,-1.5) {$v_5$};
            \node (v6) at (2,-2.5) {$v_6$};

            \draw (u1) -- (v2);
            \draw (u1) -- (v4);
            \draw (u3) -- (v4);
            \draw (u3) -- (v6);
            \draw (u5) -- (v6);
        \end{tikzpicture}
        \caption{\scriptsize Bipartite Conversion}
    \end{subfigure}
    \hfill
    \begin{subfigure}[b]{0.23\textwidth}
        \centering
        \begin{tikzpicture}[->,>=stealth, node distance=0.7cm and 0.5cm]
            \tikzstyle{every node} = [draw, circle, minimum size=5mm, font=\scriptsize]
            
            \node (u1) at (0,2.5) {$u_1$};
            \node (u2) at (0,1.5) {$u_2$};
            \node (u3) at (0,0.5) {$u_3$};
            \node (u4) at (0,-0.5) {$u_4$};
            \node (u5) at (0,-1.5) {$u_5$};
            \node (u6) at (0,-2.5) {$u_6$};

            \node (v1) at (2,2.5) {$v_1$};
            \node (v2) at (2,1.5) {$v_2$};
            \node (v3) at (2,0.5) {$v_3$};
            \node (v4) at (2,-0.5) {$v_4$};
            \node (v5) at (2,-1.5) {$v_5$};
            \node (v6) at (2,-2.5) {$v_6$};

            \draw (u1) -- (v2);
            \draw (u3) -- (v6);
        \end{tikzpicture}
        \caption{\scriptsize Maximum matching}
    \end{subfigure}
    \hfill
    \begin{subfigure}[b]{0.23\textwidth}
        \centering
        \begin{tikzpicture}[->,>=stealth, node distance=1cm and 1cm]
            \tikzstyle{every node} = [draw, circle, minimum size=5mm, font=\scriptsize]

            \node (1) at (0,0) {1};
            \node (2) [right=of 1] {2};
            \node (3) [below=of 1] {3};
            \node (4) [right=of 3] {4};
            \node (5) [below=of 3] {5};
            \node (6) [right=of 5] {6};

            \draw[very thick, red] (1) -- (2);
            \draw[very thick, blue] (3) -- (6);
        \end{tikzpicture}
        \caption{\scriptsize Minimum Path Cover}
    \end{subfigure}

    \caption{Illustration of rake minimization via minimum path cover and bipartite matching}
    \label{fig:mpc_matching}
\end{figure}
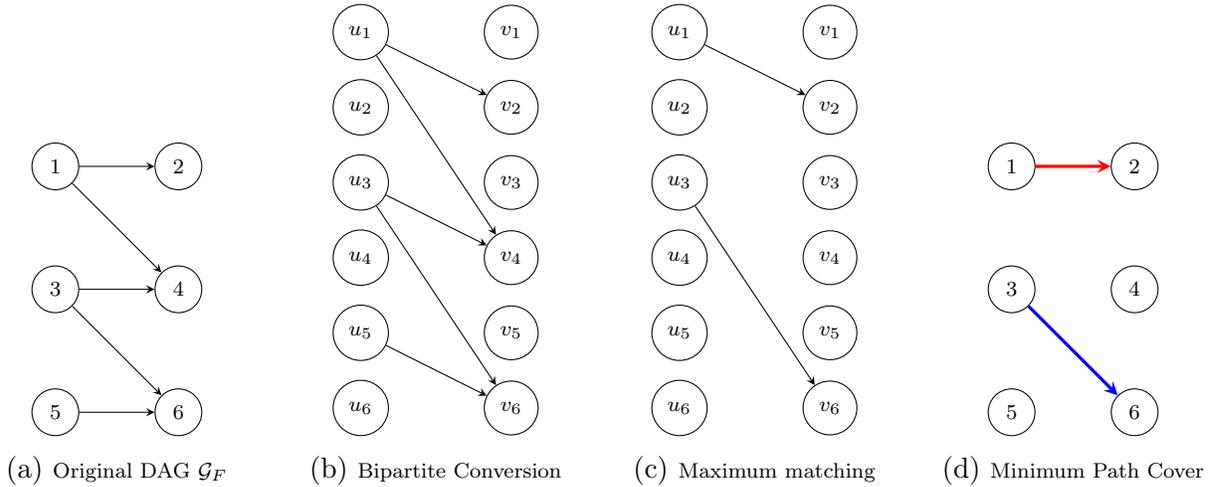

\vspace{1em}
Fig.\ref{fig:mpc_matching} illustrates the process with a simple example.

\begin{enumerate}[label=(\alph*)]
  \item \textbf{Original DAG (\(\mathcal{G}_F\))}: The nodes \(1\)-\(6\) represent individual services, and the directed edges represent the feasibility of the link.

  \item \textbf{Bipartite Conversion}:  
    Each service \(i\) is duplicated as \(u_i\) on the left partition and \(v_i\) on the right. Every feasible link \(i \to j\) in the DAG becomes the directed edge \((u_i,\,v_j)\).

  \item \textbf{Maximum Matching}:  
    The We perform Maximum Bipartite Matching and the selected edges were: \((u_1,\,v_2)\) and \((u_3,\,v_6)\).

  \item \textbf{Minimum Path Cover}:  
    Each matched pair corresponds to a continuous path in the original DAG (namely \(1 \to 2\) and \(3 \to 6\)). Unmatched vertices initiate singleton paths. Together, these four vertex-disjoint paths cover all six services, which implies that only four rakes are required.
\end{enumerate}

\subsection{Multi-Objective Variant}
\label{sec:multi_obj}

We now extend the formulation to a \textbf{multi-objective} variant where several operational objectives are considered simultaneously. In a time table it is possible for some services to be completely isolated with respect to their origin and destination; meaning no other services directly connect to either of its origin and destination stations. This means a deadhead movement is indispensable to bring a physical rake to the origin of this isolated service. A deadhead move refers to a rake traveling empty from the terminal of one service to the origin of the next service. Although such movements increase total distance and operating cost, they enable additional links and can help reduce the overall size of the fleet. Furthermore, unbalanced assignments on rakes (in terms of number of services or total distance run) can cause operational fatigue, inefficient utilization, and maintenance problems. Similarly, long waiting times between services can reduce network responsiveness and platform utilization. Keeping all these in mind, we formulate a multi-objective variant.

\subsubsection*{Problem Formulation}

Let \(\mathcal{S} = \{1,2,\dots,N\}\) be the set of scheduled services, and let \(\mathcal{R}\) be the set of rakes used in a given assignment.  Each rake \(r \in \mathcal{R}\) operates a \emph{link} 
\[
L_r = \bigl(i_{r,1},\,i_{r,2},\,\dots,\,i_{r,\,|L_r|}\bigr),
\]
where \(i_{r,k} \in \mathcal{S}\) are the services assigned to rake \(r\) in chronological order. To formulate our optimization objectives, we first define a few essential metrics associated with a rake-link \(L_r\). The \emph{length} of a link, denoted \(\lvert L_r \rvert\), refers to the total number of services assigned to rake \(r\):
\[
\lvert L_r\rvert \;=\; \text{(number of services assigned to }r\text{)}.
\]

For any two consecutive services \(i \to j\) within a link, we define the \emph{headway} \(h_{ij}\) as the time interval between the arrival of service \(i\) and the departure of service \(j\):
\[
h_{ij} = t^{\mathrm{dep}}_{j} - t^{\mathrm{arr}}_{i}.
\]

In cases where the rake must travel without passengers between the terminal of service \(i\) and the origin of service \(j\), the associated \emph{deadhead distance} is defined as:
\[
d_{ij} \;=\; \mathrm{dist}\bigl(d_i \to o_j\bigr),
\]
where \(d_i\) is the destination station of service \(i\), and \(o_j\) is the origin station of service \(j\).

Lastly, we define \emph{course length} \(c_r\) of the rake \(r\) as the total distance traveled through all assigned services, excluding any deadhead movements:
\[
c_r \;=\; \sum_{k=1}^{\lvert L_r\rvert} \mathrm{dist}\bigl(\text{origin of }i_{r,k}\;\rightarrow\;\text{destination of }i_{r,k}\bigr).
\]

These metrics serve as the building blocks for the objectives we aim to optimize in the subsequent formulation.
\clearpage

\paragraph{Objectives:}
The following are the five objectives we seek to optimize.

\begin{enumerate}
  \item \textbf{Fleet size:}
    \[
      f_{1} \;=\; \lvert \mathcal{R}\rvert.
    \]

  \item \textbf{Maximum headway across all links:}
    \[
      f_{2}
      \;=\;
      \max_{r \in \mathcal{R}}\ 
        \max_{\substack{(i,j)\in L_r}}
        \;h_{ij}.
    \]

  \item \textbf{Maximum deadhead distance across all links:}
    \[
      f_{3} = \max_{r \in \mathcal{R}}\ \max_{\substack{(i,j)\in L_r}} d_{ij}.
    \]

  \item \textbf{Link‐length standard deviation:}
    \[
      f_{4}
      \;=\;
      \mathrm{StdDev}\Bigl(\bigl\lvert L_r\bigr\rvert : r \in \mathcal{R}\Bigr).
    \]

  \item \textbf{Course‐length standard deviation:}
    \[
      f_{5}
      \;=\;
      \mathrm{StdDev}\bigl(c_r : r \in \mathcal{R}\bigr).
    \]
\end{enumerate}

Subject to the same constraints as in the single objective variant.

In this study we do not solve the multi objective optimization problem, but rather use constraint bounds to generate different link feasibility graphs and find the minimum fleet size in each of them. Then we calculate the other objectives in each case and present a comparative analysis of the results. In a future study we intend to perform a true multi objective optimization using these objectives.

\section{Results}
\label{sec:results}
We choose four key decision bounds: $w_{\min}$ as the minimum waiting time,  $w_{\max}$ as the maximum waiting time, $d_{\max}$ as the maximum deadhead distance and $v_{\mathrm{avg}}^{\max}$ as the maximum average deadhead speed. The following bounds are chosen for these parameters to act as constraints for our decision variables:
\begin{description}
  \item[$w_{\min}$] $\in \{\,0,\;60,\;120,\;180,\;240,\;300,\;\infty\}$
  \item[$w_{\max}$] $\in \{\,360,\;420,\;480,\;\dots,\;3540,\;3600,\;\infty\}$
  \item[$d_{\max}$] $\in \{\,0,\;5,\;10,\;15,\;20,\;25,\;30,\;35,\;40,\;45,\;50,\;51,\;\infty\}$
  \item[$v_{\mathrm{avg}}^{\max}$] $\in \{\,10,\;20,\;30,\;40,\;50,\;60,\;\infty\}$.
\end{description}

Here, infinity serves as a diagnostic bound to investigate theoretical limits of the model; in practice, such values lie outside the feasible region but allow us to observe boundary‐case behavior. We construct the cartesian product of all choices for $\bigl(w_{\min},\,w_{\max},\,d_{\max},\,v_{\mathrm{avg}}^{\max}\bigr)$
This yields a raw set of parameter tuples. We then remove any tuple for which $w_{\max} <= w_{\min}$
since the maximum allowed waiting time must never be less than the minimum waiting time. After excluding these infeasible combinations (including those involving \(\infty\) bounds that violate the above inequality), the resulting decision‐space comprises 44352 distinct parameter combinations that decide the design of the link feasibility graph. We then find the minimum path cover for each combination and the following observations are presented.

\subsection{Minimum theoretical fleet size}
\begin{figure}[htbp]
    \centering
    \includegraphics[width=1\linewidth]{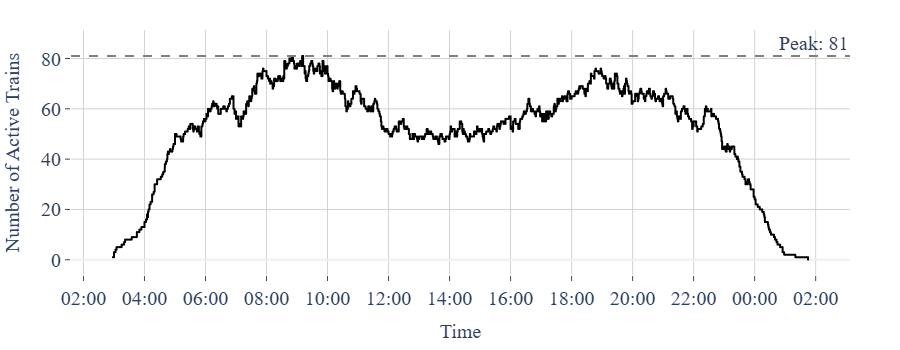}
    \caption{Number of rakes live during every 1-second duration throughout the entire day}
    \label{fig:rake_density_per_second}
\end{figure}

Throughout the day, different rakes execute different services and if we count the number of services that are live at any particular duration of the day, the maximum of this value will give us the lower limit of the theoretical bare minimum number of physical rakes necessary to run all services, since in that duration each service must be catered for by a distinct rake. Fig.\ref{fig:rake_density_per_second} shows the active rake density per second throughout the entire day. The characteristic bumps in morning and evening rush hours are distinctly visible in this plot.

The highest number of services that are live simultaneously during any time of the day is 81 (corresponding to the highest value in the y-axis of Fig.\ref{fig:rake_density_per_second}). In our experiments a minimum fleet size of 82 was achieved in 24 out of the 44352 constraint combinations discussed earlier. This corroborates our assertion that the theoretical minimum achieved through the minimum path cover, 82, must be lower bounded by the highest live rake density per second, which 81. All parameter combinations leading to fleet size of 82 and their corresponding objective values are shown in Table \ref{tab: minimum_theoretical_fleet}. In all of these cases, the $v^{max}_{avg}$ is $\infty$ indicating that in order to achieve this minimum fleet size, some rakes might need to instantly travel to the start station of the next service with infinite speed, which is absurd. 

\vspace{1em}

{\footnotesize
\begin{longtable}{|rrrrrrrrr|}

\hline
$w_{min}$   & $w_{max}$   & $d_{max}$ & $v^{max}_{avg}$ & $f_1$   & $f_2$   & $f_3$   & $f_4$ & $f_5$ \\
\hline
\endfirsthead

\multicolumn{9}{c}%
  {{\bfseries \tablename\ \thetable{} -- continued from previous page}} \\
\hline
$w_{min}$   & $w_{max}$   & $d_{max}$ & $v^{max}_{avg}$ & $f_1$    & $f_2$   & $f_3$   & $f_4$ & $f_5$ \\
\hline
\endhead

\hline \multicolumn{9}{r}{{Continued on next page}} \\ \hline
\endfoot

\endlastfoot

 0.0  & 3240.0 & inf  & inf  & \cellcolor{yellow}82 & 19500.0 & 387.660 & 2.837701 & 119.666217 \\
 0.0  & 3300.0 & inf  & inf  & \cellcolor{yellow}82 & 22560.0 & 421.476 & 2.529471 & 112.464119 \\
 0.0  & 3360.0 & inf  & inf  & \cellcolor{yellow}82 & 21960.0 & 403.280 & 2.667259 & 111.177215 \\
 0.0  & 3420.0 & inf  & inf  & \cellcolor{yellow}82 & \cellcolor{yellow}18780.0 & 452.406 & \cellcolor{yellow}2.445091 & 112.705799 \\
 0.0  & 3480.0 & inf  & inf  & \cellcolor{yellow}82 & 20340.0 & 459.426 & 2.615849 & 116.155823 \\
 0.0  & 3540.0 & inf  & inf  & \cellcolor{yellow}82 & 19680.0 & 470.690 & 2.480183 & 118.183048 \\
 0.0  & 3600.0 & inf  & inf  & \cellcolor{yellow}82 & 18960.0 & 615.020 & 2.587377 & 114.799649 \\
 0.0  &  inf   & inf  & inf  & \cellcolor{yellow}82 & 31380.0 & \cellcolor{yellow}287.730 & 2.455169 & 109.186139 \\
60.0  & 3240.0 & inf  & inf  & \cellcolor{yellow}82 & 19500.0 & 387.660 & 2.837701 & 119.666217 \\
60.0  & 3300.0 & inf  & inf  & \cellcolor{yellow}82 & 22560.0 & 421.476 & 2.529471 & 112.464119 \\
60.0  & 3360.0 & inf  & inf  & \cellcolor{yellow}82 & 21960.0 & 403.280 & 2.667259 & 111.177215 \\
60.0  & 3420.0 & inf  & inf  & \cellcolor{yellow}82 & \cellcolor{yellow}18780.0 & 452.406 & \cellcolor{yellow}2.445091 & 112.705799 \\
60.0  & 3480.0 & inf  & inf  & \cellcolor{yellow}82 & 20340.0 & 459.426 & 2.615849 & 116.155823 \\
60.0  & 3540.0 & inf  & inf  & \cellcolor{yellow}82 & 19680.0 & 470.690 & 2.480183 & 118.183048 \\
60.0  & 3600.0 & inf  & inf  & \cellcolor{yellow}82 & 18960.0 & 615.020 & 2.587377 & 114.799649 \\
60.0  &  inf   & inf  & inf  & \cellcolor{yellow}82 & 31380.0 & \cellcolor{yellow}287.730 & 2.455169 & 109.186139 \\
120.0 & 3240.0 & inf  & inf  & \cellcolor{yellow}82 & 19500.0 & 456.820 & 2.785004 & 117.922091 \\
120.0 & 3300.0 & inf  & inf  & \cellcolor{yellow}82 & 21960.0 & 421.516 & 2.657986 & 119.506923 \\
120.0 & 3360.0 & inf  & inf  & \cellcolor{yellow}82 & 21960.0 & 418.040 & 2.798271 & 121.623388 \\
120.0 & 3420.0 & inf  & inf  & \cellcolor{yellow}82 & 18840.0 & 421.516 & 2.671884 & 119.062872 \\
120.0 & 3480.0 & inf  & inf  & \cellcolor{yellow}82 & 21900.0 & 454.800 & 2.713151 & 118.914571 \\
120.0 & 3540.0 & inf  & inf  & \cellcolor{yellow}82 & 21900.0 & 403.280 & 2.859371 & 121.910674 \\
120.0 & 3600.0 & inf  & inf  & \cellcolor{yellow}82 & 20520.0 & 615.020 & 2.671884 & 119.163204 \\
120.0 &  inf   & inf  & inf  & \cellcolor{yellow}82 & 29340.0 & \cellcolor{yellow}287.730 & 2.548919 & \cellcolor{yellow}98.175734  \\
\hline
\caption{All parameters combinations for which fleet size of 82 was achieved and corresponding objective values. For each objective the cell with minimum value is highlighted.}
\label{tab: minimum_theoretical_fleet}\\
\end{longtable}
}

\subsection{Minimum achievable fleet size}

In our experiments with various combinations of parameters for the decision bounds, lower values were achieved for each objective or combinations thereof. Since there are five objectives, there are 31 possible combinations of the objectives (excluding the null combination).

\begin{table}[ht]
\centering
\footnotesize
\begin{tabular}{ll |ll |ll}
\hline
\textbf{Combo} & \textbf{Count} & \textbf{Combo} & \textbf{Count} & \textbf{Combo} & \textbf{Count} \\
\hline
($f_2$) & 43939 & ($f_3$) & 33755 & ($f_2$, $f_3$) & 33530 \\
($f_4$) & 26125 & ($f_2$, $f_4$) & 25712 & ($f_3$, $f_4$) & 23286 \\
($f_2$, $f_3$, $f_4$) & 23061 & ($f_5$) & 14641 & ($f_4$, $f_5$) & 14641 \\
($f_2$, $f_5$) & 14476 & ($f_2$, $f_4$, $f_5$) & 14476 & ($f_3$, $f_5$) & 13854 \\
($f_3$, $f_4$, $f_5$) & 13854 & ($f_2$, $f_3$, $f_5$) & 13831 & ($f_2$, $f_3$, $f_4$, $f_5$) & 13831 \\
($f_1$) & 817 & ($f_1$, $f_2$) & 418 & ($f_1$, $f_4$) & 721 \\
($f_1$, $f_5$) & 319 & ($f_1$, $f_4$, $f_5$) & 319 & ($f_1$, $f_2$, $f_4$) & 322 \\
($f_1$, $f_2$, $f_5$) & 156 & ($f_1$, $f_2$, $f_4$, $f_5$) & 156 & ($f_1$, $f_3$) & 249 \\
($f_1$, $f_3$, $f_4$) & 249 & ($f_1$, $f_2$, $f_3$) & 37 & ($f_1$, $f_2$, $f_3$, $f_4$) & 37 \\
($f_1$, $f_3$, $f_5$) & 22 & ($f_1$, $f_3$, $f_4$, $f_5$) & 22 & ($f_1$, $f_2$, $f_3$, $f_5$) & 1 \\
\multicolumn{4}{c}{} & ($f_1$, $f_2$, $f_3$, $f_4$, $f_5$) & 1 \\
\hline
\end{tabular}
\caption{Number of solutions with lower values for an objective of combinations thereof.}
\label{tab:combo_counts_rowwise}
\end{table}

We calculate all the objectives using the time table, network topology and established rake-link and the results are presented in Table \ref{tab:combo_counts_rowwise}. The significant amount of solutions in most combinations is an indication that if the network planners are to concentrate on any particular objective or combinations thereof, there are many solutions to choose from. However, the most significant achievement of these experiments is conveyed in the last row, which shows that there is one solution found which is better in all the objectives than the current existing rake-link that is used. The railway division under study uses a fleet of 99 rakes to carry out 887 suburban services throughout the day. For clarity, this solution is presented in Table \ref{tab:better_rakelink} together with the objective values associated with the existing rake-link.

\vspace{1em}
\begin{table}[h]
    \centering
    \footnotesize
    \begin{tabular}{|c|c|c|c|c|c|}\hline
        \textbf{Objective} &\textbf{$f_1$} & \textbf{$f_2$} & \textbf{$f_3$} & \textbf{$f_4$} & \textbf{$f_5$} \\ \hline
        \textbf{Existing} & 99 & 37920 & 37.44 & 3.7296 & 129.61827\\ \hline
        \rowcolor{yellow}
        \textbf{New} & 94 & 35940 & 35.95 & 3.042387 & 127.641947\\ \hline
    \end{tabular}
    \caption{Objectives from existing link as compared to the new minimum achieved}
    \label{tab:better_rakelink}
\end{table}

This represents the main contribution of this study. In particular, it demonstrates that modifying the existing rake‐link configuration yields improvements across all evaluated objectives. In particular:

\begin{itemize}
\item \textbf{Fleet‐size reduction:}
\begin{itemize}
\item The number of rakes decreases from 99 to 94 (a reduction of five units).
\item This frees up five rakes to accommodate additional services without requiring new acquisitions.
\end{itemize}

\item \textbf{Adherence to maintenance schedules:}
\begin{itemize}
\item With fewer rakes in active service, planned maintenance routines can proceed uninterrupted, reducing the risk of unscheduled downtime.
\end{itemize}

\item \textbf{Cost savings:}
\begin{itemize}
\item Each rake typically incurs several thousand dollars in procurement and lifetime maintenance costs.
\item Eliminating five redundant rakes translates into substantial savings over the operational horizon of the fleet.
\end{itemize}
\end{itemize}

\subsection{Pareto Fronts}
The preceding section demonstrated a configuration that outperforms the existing rake-link on all evaluated objectives, naturally raising the question whether there are additional such solutions. In a multi‐objective optimization context, a Pareto‐optimal solution is one whose objective values cannot be improved in any dimension without sacrificing performance in at least one other. After excluding all solutions with an unbounded (infinite) deadhead maximum average speed constraint, the remaining 40320 feasible solutions were examined, yielding nineteen distinct Pareto fronts. Table \ref{tab:pareto_summary} summarizes these objective wise minima for each Pareto front.

\begin{table}[ht]
\centering
\footnotesize
\begin{tabular}{|c|c|c|c|c|c|}
\hline
\textbf{Front}
  & \textbf{$\min f_1$}
  & \textbf{$\min f_2$}
  & \textbf{$\min f_3$}
  & \textbf{$\min f_4$}
  & \textbf{$\min f_5$} \\ \hline

1  &  92 &  720  &   0  & 0.82  &  43.75 \\
2  &  92 & 1200  &   0  & 1.13  &  61.12 \\
3  &  92 & 1440  &   0  & 1.22  &  68.07 \\
4  &  92 & 4260  &   0  & 1.67  &  84.44 \\
5  &  92 & 4440  &   0  & 1.73  &  88.13 \\
6  &  92 & 4440  &   4.53  & 1.73  &  88.13 \\
7  &  92 & 7440  &   9.46  & 2.16  & 108.09 \\
8  &  92 & 7440  &   9.46  & 2.18  & 109.50 \\
9  &  92 & 7440  &   9.46  & 2.18  & 109.56 \\
10 &  92 &11700  &   9.46  & 2.70  & 115.45 \\
11 &  95 &13080  &   9.46  & 2.77  & 118.85 \\
12 &  98 &16500  &   9.46  & 2.81  & 126.61 \\
13 &  98 &19440  &  13.80  & 2.82  & 128.42 \\
14 & 128 &19440  &  14.49  & 4.31  & 203.17 \\
15 & 134 &19440  &  18.55  & 4.31  & 206.44 \\
16 & 134 &19440  &  18.55  & 4.32  & 206.96 \\
17 & 134 &19740  &  18.55  & 4.32  & 207.56 \\
18 & 138 &19980  &  18.55  & 4.34  & 207.60 \\
19 & 138 &19980  &  19.23  & 4.34  & 208.69 \\
\hline
\end{tabular}
\caption{Objective wise minima per Pareto Front}
\label{tab:pareto_summary}
\end{table}

We can observe how the minima progressively increase as we move to higher fronts, which is a typical behavior of pareto fronts. 

\begin{figure}[htbp]
    \centering
    \includegraphics[width=0.9\linewidth]{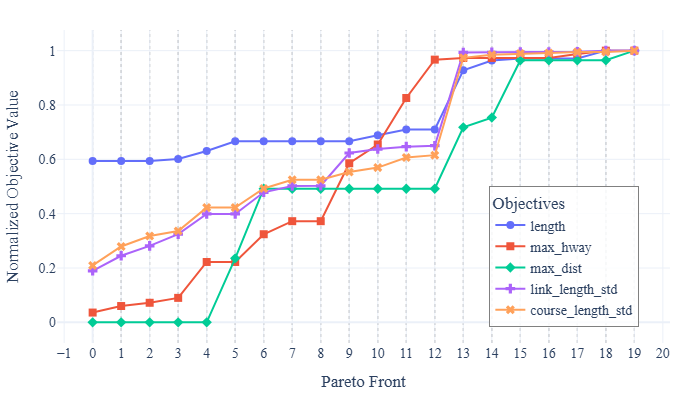}
    \caption{The minima per objective progressively increases with increasing front}
    \label{fig:pareto_minima}
\end{figure}

\vspace{1em}
The main observations from Fig.\ref{fig:pareto_minima} are:
\begin{itemize}
    \item Beyond front 12, most objectives increase sharply.
    \item max\_dist: exhibits stable intervals interrupted by sudden jumps, indicating minimal gains across certain front transitions.
    \item link\_length\_std and course\_length\_std: both rise steadily with each successive front, so even small advances markedly impact these objectives.
    \item length: remains flat through front 9, then increases modestly until front 12, jumps sharply at front 13, and stabilizes again from front 14 onward.
\end{itemize}

It may be noted that the rows in Table \ref{tab:pareto_summary} are not actual solutions, but each cell denotes the minima achieved in that particular front and in that particular objective. For example, in the first row, $f_1$ achieves minima at 92, but the corresponding other objectives can be much higher than the rest of the other values in the first row. However, the essence of a pareto front is that it is up to the end user to decide on the desired solution, depending on what factors they place their importance on: during initial investment phase, may be more emphasis is placed on the fleet size minimization, but later on, once a certain fleet has already been invested upon, the focus may shift to minimizing deadhead distance.

\subsection{Local clusters in objective space}
The Pareto fronts obtained from our multiobjective optimization reveal an important phenomenon: the presence of local clusters along each front. These clusters consist of decision configurations that, despite their variation in the decision space, correspond to nearly identical outcomes in the objective space. This observation has important implications for both theoretical understanding and practical planning. Within the 40320 solutions that were analyzed, 10333 unique objective configurations were found, which means there is a many-to-one mapping from the decision space to the objective space. Table \ref{tab:front_clusters} shows the total number of solutions and the number of clusters in the first 16 fronts.

\vspace{1em}
\begin{table}[ht]
\centering
\scriptsize
\begin{tabular}{rrrrrrrrrrrrrrrrrrrr}
\hline
\textbf{Front} & 1 & 2 & 3 & 4 & 5 & 6 & 7 & 8 & 9 & 10 & 11 & 12 & 13 & 14 & 15 & 16 \\
\hline
\textbf{count} & 15591 & 7886 & 6099 & 3771 & 2432 & 1690 & 1146 & 762 & 436 & 237 & 73 & 36 & 26 & 72 & 17 & 12 \\
\textbf{n\_clusters} & 1622 & 1360 & 1282 & 910 & 719 & 552 & 405 & 242 & 151 & 87 & 42 & 16 & 7 & 21 & 7 & 7 \\
\hline
\end{tabular}
\caption{Total number of points and number of clusters per front for the top 16 fronts}
\label{tab:front_clusters}
\end{table}

On the one hand, these clusters imply robustness in the observation space: a network controller or planner has the flexibility to adjust or perturb decision variables without adversely affecting the quality of the solution in terms of the defined objectives. This flexibility can be particularly valuable in real-world settings where decisions are subject to operational uncertainties, implementation delays, or practical constraints. In such cases, having multiple equivalent alternatives provides resilience and adaptability. For example, if we look at cluster 5, it has different values for the $v^{max}_{avg}$ but all objectives are the same. This means increasing the $v^{max}_{avg}$ from 35 to 45, won't provide any additional gains in terms of the given objectives.

On the other hand, this characteristic also reveals a potential inefficiency: redundancy in the decision space. Specifically, several combinations of decision variables may be functionally equivalent in terms of their impact on the objectives. This means that certain changes, while technically valid, do not contribute to any meaningful improvement in outcomes. For a decision-maker seeking optimal performance, such neutrality can obscure the path to actual gains and lead to wasted effort.

In this experiment, solutions with identical values in objectives have been considered in the same cluster. In practice, this can be relaxed by some $\varepsilon$-neighborhood where any solution falling within this neighborhood of each other will be a considered as part of the same cluster. This hyperparameter $\varepsilon$ will then decide the robustness of the solutions.

\begin{table}[ht]
\centering
\footnotesize
\begin{tabular}{r rrrrrrrrr}
\hline
\textbf{cluster} & \textbf{$w_{min}$} & \textbf{$w_{max}$} & \textbf{$d_{max}$} & \textbf{$v^{max}_{avg}$} & \textbf{$f_1$} & \textbf{$f_2$} & \textbf{$f_3$} & \textbf{$f_4$} & \textbf{$f_5$} \\
\hline
\rowcolor{red!10}    1 & 0.000   & 2940.000 & 20.000 & 35.000 & 135.000 & 20220.000 & 25.510 & 4.411 & 207.589 \\
\rowcolor{red!10}    1 & 60.000  & 2940.000 & 20.000 & 35.000 & 135.000 & 20220.000 & 25.510 & 4.411 & 207.589 \\
\rowcolor{blue!10}   2 & 0.000   & 3000.000 & 15.000 & 20.000 & 137.000 & 20400.000 & 18.550 & 4.356 & 208.102 \\
\rowcolor{blue!10}   2 & 60.000  & 3000.000 & 15.000 & 20.000 & 137.000 & 20400.000 & 18.550 & 4.356 & 208.102 \\
\rowcolor{green!10}  3 & 120.000 & 2820.000 & 20.000 & 30.000 & 141.000 & 19500.000 & 31.340 & 4.477 & 206.955 \\
\rowcolor{yellow!10} 4 & 240.000 & 3000.000 & 15.000 & 35.000 & 143.000 & 20640.000 & 18.550 & 4.389 & 207.609 \\
\rowcolor{yellow!10} 4 & 240.000 & 3000.000 & 15.000 & 40.000 & 143.000 & 20640.000 & 18.550 & 4.389 & 207.609 \\
\rowcolor{cyan!10}   5 & 240.000 & 3060.000 & 15.000 & 35.000 & 141.000 & 20760.000 & 18.550 & 4.370 & 207.300 \\
\rowcolor{cyan!10}   5 & 240.000 & 3060.000 & 15.000 & 40.000 & 141.000 & 20760.000 & 18.550 & 4.370 & 207.300 \\
\rowcolor{cyan!10}   5 & 240.000 & 3060.000 & 15.000 & 45.000 & 141.000 & 20760.000 & 18.550 & 4.370 & 207.300 \\
\rowcolor{orange!10} 6 & 240.000 & 3240.000 & 15.000 & 55.000 & 136.000 & 19440.000 & 19.230 & 4.318 & 207.087 \\
\rowcolor{purple!10} 7 & 240.000 & 3360.000 & 15.000 & 55.000 & 134.000 & 19440.000 & 19.230 & 4.319 & 207.857 \\
\hline
\end{tabular}
\caption{Clusters in front 16. Note that each cluster has identical objective values but different decision bounds.}
\label{tab:front16_clusters}
\end{table}

Therefore, from the perspective of planners and operational stakeholders, it is imperative to distinguish between actionable and non-impactful decisions. Understanding the structure of these local clusters enables more informed navigation of the decision space, guiding planners toward those configurations that offer tangible improvements in performance while avoiding redundant or ineffective adjustments. In this light, the identification of local clusters along the Pareto front not only provides insight into the geometry of the solution space but also becomes a practical tool for decision support in complex real-world systems.

\section{Conclusion}
\label{sec:conclusion}

This paper presents a streamlined approach to the rake-link optimization problem by applying the Minimum Path Cover (MPC) algorithm to a Directed Acyclic Graph (DAG) constructed from a topological ordering of train service departure and arrival events. Using the inherent acyclicity of the feasibility graph, the method enables computation of the minimum fleet size in polynomial time, offering a practical and scalable solution for real-time decision support.

Although the analysis incorporates multiple performance metrics, the current formulation does not constitute a true multiobjective optimization. Nevertheless, the study lays the groundwork for a future framework that will adopt a fully multi-objective paradigm that progressively converges towards Pareto-optimal fronts. A key contribution of this work is the empirical identification of one solution that strictly dominates the current solution, highlighting the potential for substantial operational improvement.

Another important finding is the presence of robust observation space clusters alongside redundant decision space clusters. This suggests that while certain decision configurations yield stable objective values, others do not offer meaningful improvement. This insight opens promising avenues for future research in decision-space sensitivity and dimensionality reduction.

The proposed DAG-based formulation and its integration with the MPC algorithm provide a computationally efficient foundation for deployment in digital twin environments, where decision alternatives can be audited and stress-tested virtually before implementation in the physical railway network. Finally, this paper expands the scope of rake-link optimization beyond traditional fleet minimization by incorporating additional operational criteria such as deadhead travel, slack time, and utilization balance, thus offering a more comprehensive and realistic approach to rail scheduling.

\bibliographystyle{unsrt}
\bibliography{main}

\end{document}